\documentclass[12pt]{amsart}

\usepackage{amsmath,xspace,amssymb,mathrsfs}
\usepackage{color}

\input xy
\xyoption{all}
\xyoption{2cell}
\UseAllTwocells
\CompileMatrices

\newcommand{\Spec}{\operatorname{Spec}}
\renewcommand{\phi}{\varphi}

\newcommand{\Ker}{\operatorname{Ker}}

\newcommand{\Min}{\operatorname{Min}}

\newcommand{\Max}{\operatorname{Max}}

\newcommand{\Ann}{\operatorname{Ann}}

\newcommand{\End}{\operatorname{End}}

\newtheorem{proposition}{Proposition}[section]
\newtheorem{lemma}[proposition]{Lemma} 
\newtheorem{corollary}[proposition]{Corollary}
\newtheorem{theorem}[proposition]{Theorem}

\newtheorem{prop-def}[proposition]{Proposition and definition}

\theoremstyle{definition}

\newtheorem{example}[proposition]{Example}

\begin{document}

\title[N-pure ideals]{Some results
on N-pure ideals}

\author[M. Aghajani]{Mohsen Aghajani}
\address{Mohsen Aghajani, Department of Mathematics, Faculty of Basic Sciences, University of Maragheh, P. O. Box 55136-553, Maragheh, Iran}
\email{aghajani14@gmail.com}


\begin{abstract}
In this paper, we consider the N-pure notion. An ideal $I$ of a ring $R$ is said to be N-pure, if for every $a\in I$ there exists $b\in I$ such that $a(1-b)\in N(R)$, where N(R) is nil radical of $R$. We provide new characterizations for N-pure ideals. In addition, N-pure ideals of an arbitrary ring are identified. Also, some other properties of N-pure ideals are studied. finally, we prove some results about the endomorphism ring of pure and N-pure ideals.\\
\end{abstract}

\maketitle

\section{\bf{Introduction}}

Throughout this paper, all rings are commutative with identity. The N-pure ideal was introduced by present author in \cite{Aghajani}. N-pure notion generalizes the pure concept. Pure ideals have been studied in various articles in the literature over the years, see e.g. \cite{Aghajani and Tarizadeh}, \cite{Al Ezeh1}, \cite{Al Ezeh2}, \cite{Borceux}, \cite{De Marco}, \cite{Fieldhouse} and \cite{Tarizadeh and Aghajani1}. In \cite{Aghajani}, it is shown that a ring is reduced if and only if its Pure and N-pure ideals are the same. Recall that a ring $R$ is said to be a mp-ring if every its prime ideal contains a unique minimal prime ideal of $R$. Using N-pure notion, it was proved that there is a new proper subclass of mp-rings which contains strictly the class of reduced mp-rings, see \cite[Theorem 4.7]{Aghajani}, \cite[Remark 4.13]{Aghajani} and \cite[Example 4.16]{Aghajani}.  Those rings were called the mid rings.\\
Also, Proposition \ref{Proposition II} specifies N-pure ideals of an arbitrary ring. This result, in combination with Theorem \ref{Theorem II}, identifies the N-pure ideals of mp-rings. In Theorem \ref{Theorem III}, we provide a new characterization for N-pure ideals in terms of a special type of ring of fractions. In the following, we show that N-purity is stable under ring extension, see Lemma \ref{Lemma II}. Also, pure ideals of a Gelfand rings are cosiderd, see Lemma \ref{Lemma III} and Proposition \ref{Proposition IV}. In addition, we show if $I$ is a finitely generated N-pure ideal of a ring $R$, then there exists $a\in I$ such that $\sqrt{I}=\sqrt{Ra}$, see Proposition \ref{Proposition VI}. Also, it is proved that the endomorphism ring of an arbitrary power of a pure ideal is commutative, see Theorem \ref{Theorem V}. Finally, we prove a result about commutativity of endomorphism ring of some power of principal N-pure ideals. \\

\section{\bf{Preliminaries}}

Here we mention some objects which is needed in the following. An ideal $I$ of a ring $R$ is said to be a pure ideal if the canonical ring homomorphism $\pi: R\longrightarrow R/I$ is a flat homomorphism. Then $I$ is a pure ideal if and only if $\Ann(a) + I = R$ for all $a\in I$, if and only if for each $a\in I$ there exists some $b\in I$ such that $a(1 - b) = 0$.
An ideal $I$ of a ring $R$ is said to be a N-pure ideal if for each $a\in I$ there exists $b\in I$ such that $a(1-b)\in N(R)$ where $N(R)$ is the nil-radical of $R$, or equivalently for each $a\in I$ there exist $n\geqslant1$ and $b\in I$ such that $a^{n}(1-b)=0$. The following result characterizes the N-pure ideals.

\begin{theorem}\label{Theorem V} Let $R$ be a ring and $I$ be an ideal of $R$. Then the following conditions are equivalent. \\
$\mathbf{(i)}$ $I$ is a N-pure ideal. \\
$\mathbf{(ii)}$ For every $a \in I$ there exists $t\geqslant1$ such that $\Ann(a^{t})+I=R$. \\
$\mathbf{(iii)}$ $\sqrt{I}=\{a\in R|\exists n \geqslant1, \Ann(a^{n})+I=R\}$.\\
$\mathbf{(iv)}$ $\sqrt{I}$ is a N-pure ideal.\\
$\mathbf{(v)}$ There exists a unique pure ideal $J$ such that $\sqrt{I}=\sqrt{J}$.
\end{theorem}

{\bf Proof.} See \cite[Theorem 2.6]{Aghajani}. $\blacksquare$\\

Recall that a ring $R$ is said to be a mid ring if for every $a\in R$, $\Ann(a)$ is a N-pure ideal. If $\mathfrak{p}\in \Spec(R)$, we denote the natural ring homomorphism $R\longrightarrow S^{-1}R$ by $\pi_{\mathfrak{p}}$ where $S=R\setminus \mathfrak{p}$. Then a ring $R$ is a mid ring if and only if $\Ker\pi_{\mathfrak{p}}$ is a primary ideal for all $\mathfrak{p}\in \Spec(R)$, see \cite[Theorem 4.6]{Aghajani}. Finally, a ring $R$ is called a Gelfand ring (pm-ring) if every prime ideal of $R$ is contained in a unique maximal ideal.

\section{\bf{Main results}}

The next result identifies N-pure ideals in an arbitrary ring.\\

\begin{proposition}\label{Proposition II} If $I$ is a N-pure ideal of a ring $R$, then  $$\sqrt{I}=\bigcap\limits_{\mathfrak{m}\in V(I)\cap \Max(R)} \sqrt{\Ker\pi_{\mathfrak{m}}}.$$ \\
\end{proposition}

{\bf Proof.} Let $a\in \sqrt{I}$. Then there exist $n\geqslant1$ and $b\in I$ such that $a^{n}(1-b)=0$. If $\mathfrak{m}$ is a maximal ideal containing $I$, then $1-b \notin \mathfrak{m}$ and so $a\in \sqrt{\Ker\pi_{\mathfrak{m}}}$. Now, let $a\in \bigcap\limits_{\mathfrak{m}\in V(I)\cap \Max(R)} \sqrt{\Ker\pi_{\mathfrak{m}}}$. It suffices to show that $\Ann(a^{m})+I=R$ for some $m\geqslant1$, see \cite[Theorem 2.6(iv)]{Aghajani}. Assume $\Ann(a^{n})+I\neq R$ for every $n\geqslant1$. Then $J=\sum\limits_{n\geqslant1} \big(\Ann(a^{n})+I\big)$ is a proper ideal. Let $\mathfrak{n}$ be a maximal ideal which contains $J$. Hence $a\notin \sqrt{\Ker\pi_{\mathfrak{n}}}$ which is a contradiction. $\blacksquare$\\

\begin{corollary}\label{Corollary I} Let $R$ be a ring and $I$ be a N-pure ideal of $R$. If $\mathfrak{m}$ is a maximal ideal of $R$ containing $I$, then $I\subseteq \sqrt{\Ker\pi_{\mathfrak{m}}}$. $\blacksquare$\\
\end{corollary}

\begin{lemma}\label{Lemma V} Let $R$ be a ring. If $I$ is a proper N-pure ideal of $R$ containing a prime ideal $\mathfrak{p}$, then $I=\mathfrak{p}$. Especially, $I$ is a minimal prime ideal of $R$.\\
\end{lemma}

{\bf Proof.} Let $\mathfrak{p}\neq I$. Then there exists $a\in I\setminus \mathfrak{p}$. Since $I$ is N-pure, then there are $b\in I$ and $n\geqslant1$ such that $a^{n}(1-b)=0$. Hence $1-b\in \mathfrak{p}$ and so $1\in I$. This is a contradiction. Thus we have $I=\mathfrak{p}$. Additionally, this implies that $I$ is a minimal prime ideal of $R$. $\blacksquare$\\

\begin{corollary}\label{Corollary IV} The radical of a N-pure ideal $I$ is the intersection of all
minimal prime ideals of $R$ containing $I$. $\blacksquare$\\
\end{corollary}
Let $\mathfrak{p}\in\Spec(R)$. By $\Lambda(\mathfrak{p})$ we mean the set of all $\mathfrak{q}\in\Spec(R)$ such that $\mathfrak{q}\subseteq\mathfrak{p}$. \\
\begin{theorem}\label{Theorem II} Let $R$ be a mp-ring and $I$ be a proper ideal of $R$. Then $$J=\bigcap
\limits_{\mathfrak{m}\in V(I)\cap \Max(R)} \sqrt{\Ker\pi_{\mathfrak{m}}}=\bigcap\limits_{\substack{\mathfrak{m}\in V(I)\cap \Max(R) \\ \mathfrak{p}\in \Lambda(\mathfrak{m})\cap \Min(R) }} \mathfrak{p}$$ is a N-pure ideal of $R$.\\
\end{theorem}

{\bf Proof.} Let $a\in J$. Then for every maximal ideal $\mathfrak{m}$ containing $I$ there are $b_{\mathfrak{m}}\in R\setminus \mathfrak{m}$ and $t_{\mathfrak{m}}\geqslant1$ such that $a^{t_{\mathfrak{m}}}b_{\mathfrak{m}}=0$. If $\mathfrak{m}$ is a maximal ideal that does not contain $I$, then let $c_{\mathfrak{m}}\in I\setminus \mathfrak{m}$. Thus the ideal $K=\big(b_{\mathfrak{m}}, c_{\mathfrak{m}} |\mathfrak{m}\in \Max(R)\big)$ is the whole ring. Hence $1=\sum\limits_{i=1}^{m}r_{i}b_{i}+\sum\limits_{j=1}^{n}s_{j}c_{j}$ where $b_{i}=b_{\mathfrak{m_{i}}}\in R\setminus \mathfrak{m_{i}}$, $c_{j}=c_{\mathfrak{m_{j}}}\in I\setminus \mathfrak{m_{j}}$ and $r_{i},s_{j}\in R$. Now, let $t=min\{t_{i}|i=1,\cdots,m\}$ where $t_{i}=t_{\mathfrak{m_{i}}}$. Then we have $a^{t}(1-c)=0$ where $c=\sum\limits_{j=1}^{n}s_{j}c_{j}$. If $\mathfrak{m}$ is a maximal ideal containing $I$, then $1-c\notin \mathfrak{m}$ and so $a^{t}\in \Ker\pi_{\mathfrak{m}}$. On the other hand, by \cite[Theorem 6.2(ix)]{Aghajani and Tarizadeh}, there exists $u\geqslant1$ such that $\Ann(a^{tu})+\Ann((1-c)^{u})=R$. Then $e+f=1$ for $e \in \Ann(a^{tu})$ and $f \in \Ann((1-c)^{u})$. Therefore, $a^{tu}(1-f)=0$ and $f\in \Ker\pi_{\mathfrak{m}}$. This yields that $J$ is a N-pure ideal. The second equality holds by \cite[Theorem 3.3(iii)]{Aghajani}. $\blacksquare$\\

\begin{proposition}\label{Proposition III} Let $R$ be a ring and $I$ be an ideal of $R$. Then the following statements are equivalent.\\
$\mathbf{(i)}$ $I$ is a N-pure ideal.\\
$\mathbf{(ii)}$ For each ideal $J$ with $I\subseteq J$, $J$ is a N-pure ideal if and only if  $J/I$ is N-pure.\\
\end{proposition}

{\bf Proof.} $\mathbf{(i)} \Rightarrow \mathbf{(ii)} :$ If $J$ is a N-pure ideal of $R$ containing $I$, then it is clear to see that $J/I$ is a N-pure ideal of $R/I$. Conversely, let $a\in J$. Then there exist $n\geqslant1$ and $b+I\in J/I$ such that $(a+I)^{n}(1+I-(b+I))=0$. Thus $a^{n}(1-b) \in I$ and so there are $m\geqslant1$ and $c\in I$ such that $a^{mn}(1-b)^{m}(1-c)=0$. Hence $a^{mn}(1-bd)(1-c)=0$ for some $d\in R$. So $a^{mn}(1-bd-c+bcd)=0$. Therefore $J$ is a N-pure ideal, since $bd+c-bcd\in J$.\\
$\mathbf{(ii)} \Rightarrow \mathbf{(i)} :$ Since the zero ideal is a N-pure ideal, then the assertion is obtained by setting $J=I$. $\blacksquare$\\

Let $I$ be an ideal of a ring $R$. Setting $S=1+I$. Then $S$ is a multiplicative closed subset of $R$. In the following, we prove some results about relation between the ring $S^{-1}R$ and pure and N-pure ideals. The following result provides a new characterization for N-pure ideals. \\

\begin{theorem}\label{Theorem III} Let $R$ be a ring, $I$ be an ideal of $R$ and $S=1+I$. Then the following statements are equivalent.\\
$\mathbf{(i)}$ $I$ is a N-pure ideal.\\
$\mathbf{(ii)}$ $\Spec(S^{-1}R)=\{S^{-1}\mathfrak{p}| I\subseteq\mathfrak{p}\}$.\\
$\mathbf{(iii)}$ $N(S^{-1}R)=S^{-1}\sqrt{I}$.\\
$\mathbf{(iv)}$ $\sqrt{I}=\sqrt{\Ker\varphi}$ where $\varphi: R\longrightarrow S^{-1}R$ is the natural ring homomorphism.\\
$\mathbf{(v)}$ $S^{-1}R/N(S^{-1}R)\cong R/\sqrt{I}$ as $R$-algebras.\\
\end{theorem}

{\bf Proof.} $\mathbf{(i)} \Rightarrow \mathbf{(ii)} :$ Let $S^{-1}\mathfrak{p}\in \Spec(S^{-1}R)$. Then we have $\mathfrak{p}\cap S=\emptyset$. We claim that $I\subseteq \mathfrak{p}$. Otherwise, there exists $a\in I\setminus \mathfrak{p}$. Thus there are $n\geqslant1$ and $b\in I$ such that $a^{n}(1-b)=0$, because $I$ is N-pure. So $1-b\in \mathfrak{p}$ which is a contradiction. On the other hand, if $I\subseteq \mathfrak{p}$, then it is obvious that $S^{-1}\mathfrak{p}\in \Spec(S^{-1}R)$.\\ $\mathbf{(ii)} \Rightarrow \mathbf{(i)} :$ It suffices to show that for each $a\in I$ there exists $n\geqslant1$ such that $\Ann(a^{n})+I=R$, see \cite[Theorem 2.6(iii)]{Aghajani}. If not then the ideal $K=\sum\limits_{m\geqslant1}\big(\Ann(a^{m})+I\big)$ is a proper ideal. Let $\mathfrak{m}$ be a maximal ideal containing $K$. This means that $a^{m}/1\neq0$ for all $m\geqslant1$ and so $a/1\notin N(S^{-1}R)$. Then there exists $S^{-1}\mathfrak{p}\in \Min(S^{-1}R)$ such that $a/1\notin S^{-1}\mathfrak{p}$. This is a contradiction because $I\subseteq\mathfrak{p}$. Hence $I$ is a N-pure ideal.\\
$\mathbf{(i)} \Rightarrow \mathbf{(iii)} :$ In general, we have $N(S^{-1}R)\subseteq S^{-1}\sqrt{I}$. Let $a/s\in S^{-1}\sqrt{I}$. Then there exist $b\in \sqrt{I}$ and $t\in S$ such that $a/s=b/t$. Hence there are $n\geqslant1$ and $c\in I$ such that $b^{n}(1-c)=0$. So $a^{n}/s^{n}=b^{n}/t^{n}=0$. This implies $a/s\in N(S^{-1}R)$. Therefore, $N(S^{-1}R)=S^{-1}\sqrt{I}$.\\
$\mathbf{(iii)} \Rightarrow \mathbf{(i)} :$ Let $a\in I$. Then $a/1\in N(S^{-1}R)$ and so there exists $n\geqslant1$ such that $a^{n}/1=0$. Hence there exists $b\in I$ such that $a^{n}(1+b)=0$. Thus $I$ is a N-pure ideal.\\
$\mathbf{(i)} \Rightarrow \mathbf{(iv)} :$ In general, we have $\Ker\varphi\subseteq I$. Let $a\in I$. Then there exist $n\geqslant1$ and $b\in I$ such that $a^{n}(1-b)=0$. Hence, $a \in \sqrt{\Ker\varphi}$. So $\sqrt{I}=\sqrt{\Ker\varphi}$.\\
$\mathbf{(iv)} \Rightarrow \mathbf{(i)} :$ It is obvious.\\
$\mathbf{(iv)} \Rightarrow \mathbf{(v)} :$ Define $\psi : S^{-1}R\longrightarrow R/\sqrt{I}$ with $\psi(a/s)=a+\sqrt{I}$. Then this map is an algebra homomorphism. Let $a/s \in \Ker\psi$. Thus $a\in \sqrt{I}$ and so by hypothesis $a\in \sqrt{\Ker\varphi}$. So there exist $n\geqslant1$ and $b\in I$ such that $a^{n}(1+b)=0$. Hence $a/s\in N(S^{-1}R)$. Therefore, $\Ker\psi=N(S^{-1}R)$. This implies that $S^{-1}R/N(S^{-1}R)\cong R/\sqrt{I}$.\\
$\mathbf{(v)} \Rightarrow \mathbf{(i)} :$ By \cite[Theorem 2.6]{Aghajani}, it suffices to show that $\sqrt{I}$ is a N-pure ideal. Let $a\in \sqrt{I}$. We claim that there is $n\geqslant1$ such that $\Ann(a^{n})+\sqrt{I}=R$. Otherwise, the ideal $\sum\limits_{m\geqslant1}\big(\Ann(a^{m})+\sqrt{I}\big)$ is a proper ideal. Hence it contained in a maximal ideal $\mathfrak{m}$. Thus $\Ann(a^{m})\subseteq\mathfrak{m}$ and so $a/1\notin N(S^{-1}R)$. This is a contradiction. Therefore, $I$ is a N-pure ideal. $\blacksquare$\\

\begin{corollary}\label{Corollary II} Let $R$ be a ring, $I$ and $J$ be two ideals of $R$ with $S=1+I$ and $T=1+J$. If $\sqrt{I}=\sqrt{J}$ and $I$ is a N-pure ideal, then $S^{-1}R/S^{-1}\sqrt{I}\cong T^{-1}R/T^{-1}\sqrt{J}$ as algebras. $\blacksquare$\\
\end{corollary}

In the following result, the criterion (iii) is new. \\
\begin{theorem}\label{Theorem IV} Let $R$ be a ring, $I$ be an ideal of $R$ and $S=1+I$. Then the following statements are equivalent.\\
$\mathbf{(i)}$ $I$ is a pure ideal.\\
$\mathbf{(ii)}$ $S^{-1}R\cong R/I$ as $R$-algebras.\\
$\mathbf{(iii)}$ $S^{-1}I$ is the zero ideal of $S^{-1}R$.\\
\end{theorem}

{\bf Proof.} $\mathbf{(i)} \Rightarrow \mathbf{(ii)} :$ Define $\varphi : S^{-1}R\longrightarrow R/I$ with $\varphi(a/s)=a+I$. Since $I$ is a pure ideal, then $\varphi$ is a well-define algebra homomorphism. Now, let $\varphi(a/s)=\varphi(b/t)$. Then we have $a+I=b+I$ and so $ta+I=sb+I$. Thus $ta-sb\in I$. So there exists $c\in I$ such that $(ta-sb)(1-c)=0$. Then $a/s=b/t$. This yields that $\varphi$ is a monomorphism. Therefore, $\varphi$ is an isomorphism.\\
$\mathbf{(ii)} \Rightarrow \mathbf{(i)} :$ Let $a\in I$. It suffices to show that $\Ann(a)+I=R$. Otherwise, there exists a maximal ideal $\mathfrak{m}$ such that $\Ann(a)+I\subseteq\mathfrak{m}$. Then by hypothesis, $a/1\neq0$. This is a contradiction. Hence $I$ is a pure ideal.\\
$\mathbf{(i)} \Rightarrow \mathbf{(iii)} :$ Let $a/s \in S^{-1}I$. Then $a/s=b/t$ for some $b\in I$ and $t\in S$. Thus there exists $c\in I$ such that $b(1-c)=0$. Hence $a/s=0$. This yields that $S^{-1}I$ is the zero ideal of $S^{-1}R$.\\
$\mathbf{(iii)} \Rightarrow \mathbf{(i)} :$ It is clear.
$\blacksquare$\\

\begin{example}\label{example I} It is well-Known that only pure(=N-pure) ideals of an integral domain are zero ideal and the whole ring. Let $p$ be a prime number. Then $p^{n} \mathbb{Z}$ is not a Pure ideal of the ring $ \mathbb{Z}$ for $n\geqslant1$. If $S=1+p^{n} \mathbb{Z}$, then $S^{-1}\mathbb{Z}$ is an integral domain. If $q$ is a prime number with $p\neq q$, then we have $q\mathbb{Z}\cap S\neq \emptyset$. Hence $\Spec(S^{-1}\mathbb{Z})=\{0, S^{-1}p \mathbb{Z}\}$. In case $n=1$, $\mathbb{Z}/p\mathbb{Z}$ is a field while $S^{-1}\mathbb{Z}$ is a non-field integral domain. As another example, let $R=\mathbb{Z}/12\mathbb{Z}$,  $\mathfrak{p}=2\mathbb{Z}/12\mathbb{Z}$ and $S=1+\mathfrak{p}$. It is easy to check that $\mathfrak{q}=3\mathbb{Z}/12\mathbb{Z}$ is a pure ideal of $R$. Hence, by \cite[Theorem 3.4(iv)]{Aghajani} and \cite[Corollary 3.5(iv)]{Aghajani}, $\mathfrak{p}$ is a N-pure ideal of $R$ which is not pure. Clearly, $\Spec(S^{-1}R)=\{ S^{-1}\mathfrak{p}\}$. Then $S^{-1}R$ is not an integral domain, since $S^{-1}\mathfrak{p}\neq0$. But $R/\mathfrak{p}\cong\mathbb{Z}/2\mathbb{Z}$ which is a field.\\
\end{example}

\begin{corollary}\label{Corollary III} Let $R$ be a ring, $I$ be a pure (N-pure) ideal of $R$ and $S=1+I$. Then $R/I$ is a mp-ring (Gelfand) ring if and only if $S^{-1}R$ is a mp-ring (Gelfand) ring. \\
\end{corollary}

{\bf Proof.} It follows from Theorems \ref{Theorem III} and \ref{Theorem IV}. $\blacksquare$\\

\begin{lemma}\label{Lemma II} Let $\phi:R\longrightarrow S$ be a ring homomorphism. If $I$ is a N-pure ideal of $R$, then $IS$ is a N-pure ideal of $S$.\\
\end{lemma}

{\bf Proof.} Let $b\in IS$. Then there are $n\geqslant1$, $a_{i}\in I$ and $s_{i}\in S$ such that $b=\sum\limits_{i=1}^{n}a_{i}s_{i}$. Thus by \cite[Theorem 2.6(ii)]{Aghajani}, there exist $c\in I$ and $t\geqslant1$ such that $a_{i}^{t}(1-c)=0$ for $i=1,\cdots,n$. Then $$b^{nt}(1-c)=\sum\limits_{\sum\limits_{j=1}^{n} i_{j}=nt}a_{1}^{i_{1}}\cdots a_{n}^{i_{n}}s'_{i_{1},\cdots,i_{n}}(1-c)=0.$$
Therefore, $IS$ is a N-pure ideal of $S$. $\blacksquare$\\

In the next result, we give the new short proof for \cite[Lemma 1.2]{Al Ezeh2} in non-reduced case.\\

\begin{lemma}\label{Lemma III} Let $R$ be a Gelfand ring. Then $\Ker\pi_{\mathfrak{m}}$ is a pure ideal for every $\mathfrak{m}\in \Max(R)$.\\
\end{lemma}

{\bf Proof.} Let $a\in \Ker\pi_{\mathfrak{m}}$. Then there exists $b\in R\setminus \mathfrak{m}$ such that $ab=0$. It is suffices to show that $\Ann(a)+\Ker\pi_{\mathfrak{m}}=R$. Otherwise, there exists a maximal ideal $\mathfrak{m'}$ such that $\Ann(a)+\Ker\pi_{\mathfrak{m}}\subseteq\mathfrak{m'}$.
Hence $\mathfrak{m}\neq\mathfrak{m'}$, since $\Ann(a)\nsubseteq\mathfrak{m}$. But this is a contradiction by \cite[Theorem 4.3(ix)]{Aghajani and Tarizadeh}. $\blacksquare$\\

The following result improves \cite[Theorem 1.5]{Al Ezeh2}.\\

\begin{proposition}\label{Proposition IV} Let $R$ be a Gelfand ring and $\mathfrak{m}$ be a maximal ideal of $R$. Then $\mathfrak{m}$ is pure if and only if $\Ker\pi_{\mathfrak{m}}=\mathfrak{m}$.
\end{proposition}

{\bf Proof.} It follows immediately from Lemma \ref{Lemma III}. $\blacksquare$\\

\begin{lemma}\label{Lemma IV} Let $R$ be a ring. If $\mathfrak{p}$ is a N-pure prime ideal of $R$, then $\sqrt{\Ker\pi_{\mathfrak{p}}}= \mathfrak{p}$.\\
\end{lemma}

{\bf Proof.} It is clear that $\sqrt{\Ker\pi_{\mathfrak{p}}}\subseteq \mathfrak{p}$. Now, let $a\in \mathfrak{p}$. Then there exist $n\geqslant1$ and $b\in \mathfrak{p}$ such that $a^{n}(1-b)=0$. Thus $a^{n}\in \Ker\pi_{\mathfrak{p}}$ and so $a\in \sqrt{\Ker\pi_{\mathfrak{p}}}$. This means that $\sqrt{\Ker\pi_{\mathfrak{p}}}= \mathfrak{p}$. $\blacksquare$\\

\begin{lemma}\label{Lemma VI} Let $R$ be a ring and $I$ be an ideal of $R$ with $I\subseteq J(R)$. If $I$ is a N-pure ideal, then $I\subseteq N(R)$.\\
\end{lemma}

{\bf Proof.} Let $a\in I$. Then there exist $n\geqslant1$ and $b\in I$ such that $a^{n}(1-b)=0$. Since $b\in J(R)$, then $1-b$ is invertible and so $a^{n}=0$. Hence $I\subseteq N(R)$. $\blacksquare$\\

Recall that an ideal $I$ of a ring $R$ is said to be $\emph{strongly $\pi$-regular}$ if for every $a\in I$
there are $n\geqslant1$ and $b\in I$ such that $a^{n}=a^{n+1}b$. It is clear that every strongly $\pi$-regular is a N-pure ideal. Also, it is well known
that a ring is zero-dimensional if and only if it is a strongly $\pi$-regular ring.\\

\begin{proposition}\label{Proposition V} Let $R$ be a ring. Then $R$ is a zero-dimensional ring if and only if  the N-pure ideals and the strongly $\pi$-regular ideals are the same.\\
\end{proposition}

{\bf Proof.} It is suffices to show that every N-pure ideal is a strongly $\pi$-regular ideal. Let $I$ be a N-pure ideal and $a\in I$. Then by \cite[Theorem 3.4(iii)]{Aghajani}, there is $n\geqslant1$ and $c\in R$ such that $a^{n}=a^{n+1}c$. Hence $I$ is a strongly $\pi$-regular ideal. The converse case holds by \cite[Theorem 3.4(ii)]{Aghajani} and the fact that $R$ is a N-pure ideal. $\blacksquare$\\

\begin{proposition}\label{Proposition VI} Let $R$ be a ring and $I$ be a finitely generated N-pure ideal of $R$. Then $\sqrt{I}=\sqrt{Ra}$ for some $a\in I$.\\
\end{proposition}

{\bf Proof.} Let $\{a_{1},\cdots,a_{n}\}$ be a generating set of $I$. Then there exist $t\geqslant1$ and $a\in I$ such that $a_{i}^{t}(1-a)=0$ for all $i=1,\cdots,n$ by \cite[Theorem 2.6(ii)]{Aghajani}. Thus for all $i=1,\cdots,n$ we have $a_{i}^{t}=a_{i}^{t}a$. Hence $a_{i}\in \sqrt{Ra}$ for each $1\leqslant i\leqslant n$. Therefore $I\subseteq\sqrt{Ra}$ and so $\sqrt{I}=\sqrt{Ra}$. $\blacksquare$\\

\begin{lemma}\label{Lemma VII} Let $R$ be a ring and $I$ be a finitely generated N-pure ideal of $R$. Then there exist  $n\geqslant1$ and $b\in I$ such that $a^{n}(1-b)=0$ for all $a\in I$.\\
\end{lemma}

{\bf Proof.} Let $I=\sum\limits_{i=1}^{n}Ra_{i}$. Then similar to the proof of \ref{Proposition VI} there exist $t\geqslant1$ and $b\in I$ such that $a_{i}^{t}(1-b)=0$ for all $i=1,\cdots,n$. Thus if $a=\sum\limits_{i=1}^{n}r_{i}a_{i}$, then $$a^{nt}=\sum\limits_{i_{1}+\cdots+i_{n}=nt}s_{i_{1},\cdots,i_{n}}a_{1}^{i_{1}}\cdots a_{n}^{i_{n}}=
\sum\limits_{i_{1}+\cdots+i_{n}=nt}s_{i_{1},\cdots,i_{n}}a_{1}^{i_{1}}\cdots a_{n}^{i_{n}}b=a^{nt}b$$ since at least $i_{j}\geqslant t$ for some $j$. Hence the assertion is proved. $\blacksquare$\\

It is well known that the endomorphism ring of an arbitrary module is a noncommutative ring. In the following result, we prove that the endomorphism ring of arbitrary power of a pure ideal is commutative.\\

\begin{theorem}\label{Theorem V} Let $R$ be a ring and $I$ be a pure ideal of $R$. Then $\End(I^{n})$ is a commutative ring for all $n\geqslant1$.\\
\end{theorem}

{\bf Proof.} Let $n$ be a fixed positive integer and $f,g\in \End(I^{n})$. If $a_{1},\cdots,a_{n} \in I$, then there exist $b_{i}\in I$ such that $a_{i}=a_{i}b_{i}$ for all $i=1,\cdots,n$. Setting $1-b=\prod\limits_{i=1}^{n}(1-b_{i})$, then we have $b\in I$ and so $a_{i}=a_{i}b$ for all $i$. Hence
\begin{eqnarray}
f g(a_{1}\cdots a_{n})&=&f g(a_{1}\cdots a_{n}b^{n})  \nonumber \\ &=&f( g(a_{1}\cdots a_{n}b^{n}))  \nonumber \\  &=&f(a_{1}\cdots a_{n}g(b^{n})) \nonumber \\  &=&g(b^{n})f(a_{1}\cdots a_{n})  \nonumber \\  &=&g(f(a_{1}\cdots a_{n}b^{n}) \nonumber \\  &=&g f(a_{1}\cdots a_{n}b^{n}) \nonumber \\  &=&g f(a_{1}\cdots a_{n}). \nonumber
\end{eqnarray}
Then we have $fg=gf$, since $f$ and $g$ are stable under sum operation. Therefore $\End(I^{n})$ is a commutative ring. $\blacksquare$\\

\begin{lemma}\label{Lemma VII} Let $R$ be a ring and $I$ be a principal N-pure ideal of $R$. Then there exists $n\geqslant1$ such that $\End(I^{i})$ is a commutative ring for each $i\geqslant n$.\\
\end{lemma}

{\bf Proof.} Let $I=Ra$. Then there exist $n\geqslant1$ and $r\in R$ such that $a^{n}(1-ra)=0$. It is obvious that
$a^{i}=a^{i}(ra)^{i}$ for each $i\geqslant n$. Now let $i$ be a fixed integer with $i\geqslant n$ and $f,g\in \End(I^{i})$. Then we have
\begin{eqnarray}
f g(a^{i})&=&f g(a^{i}(ra)^{i})  \nonumber \\ &=&f( g(a^{i}(ra)^{i}))  \nonumber \\  &=&f(a^{i}g((ra)^{i})) \nonumber \\  &=&g((ra)^{i})f(a^{i})  \nonumber \\  &=&g(f(a^{i}(ra)^{i})) \nonumber \\  &=&g f(a^{i}(ra)^{i}) \nonumber \\  &=&g f(a^{i}) \nonumber
\end{eqnarray}
Hence $fg=gf$ and so $\End(I^{i})$ is a commutative ring for each $i\geqslant n$. $\blacksquare$\\

\end{document}